\documentclass[12pt,twoside]{amsart}
\usepackage{amssymb,amsmath}

\title{Formal neighborhoods in arc spaces}
\author{Peter Petrov}

\newcommand{\sing}{\operatorname{Sing}}
\newcommand{\spec}{\operatorname{Spec}}

\newcommand{\ord}{\operatorname{ord}}
\newcommand{\Hom}{\operatorname{Hom}}

\newcommand{\im}{\operatorname{Im}}

\newcommand{\spf}{\operatorname{Spf}}
\newcommand{\ch}{\operatorname{char}}
\newcommand{\cod}{\operatorname{codim}}

\newtheorem{thm}{Theorem}[section]
\newtheorem{lem}[thm]{Lemma}

\newtheorem{defn}[thm]{Definition}
\newtheorem{Cor}[thm]{Corolary}
\newtheorem{Prop}[thm]{Proposition}
\theoremstyle{definition}
\newtheorem{Ex}[thm]{Example}
\newtheorem{Rem}[thm]{Remark}

\begin{document}

\maketitle

\begin{abstract}

The theorem of Greenberg-Kazhdan-Drinfeld describes the formal neighborhood of a closed arc. After giving a complete proof with examples, two possible versions for the relative case of the theorem are discussed. Each one is shown to hold for some classes of morphisms, including smooth morphisms and closed embeddings in smooth variety.

\end{abstract}



\section{Introduction}
The arc spaces have been introduced by J. Nash (\cite{Nash}) in an attempt to study the singularities of an algebraic variety $X$ by attaching a geometric space intrinsically to $X$. This is the space of arcs $X_{\infty}$, which reflects the geometry of the singular locus without need to take resolution of singularities. It was introduced also in an attempt to characterize the essential divisors over $X$, i.e. those divisors appearing birationally as irreducible components of the exceptional locus in any resolution of singularities of $X$. After showing that there is an injection between the components in $X_{\infty}$ of the singular fiber, and the essential divisors over $X$, Nash proposed what is called now the Nash problem on arc spaces. The geometry of these spaces has been studied actively since then (\cite{IK}, \cite{Ishii}, \cite{ELM}, \cite{math.AG/0102201}, \cite{MR1856396}, \cite{MR1664700}, \cite{MR1359701},   \cite{MR2197405},
  \cite{MR2900541}, \cite{MR2979864}, \cite{MR2141063}, \cite{LR},
   \cite{MR2245998}, \cite{MR2604083},
   \cite{MR2359548}, H.Mourtada - Jet schemes of rational double point surface singularities, Valuation Theory in Interaction, EMS Ser. Congr. Rep., Eur. Math. Soc., 2014,
   \cite{MR2000468}, \cite{perez-2007}, \cite{MR1053487}, \cite{MR2675049}, \cite{MR3089196}, M. Leyton-Alvarez - Resolution de probleme des arcs de Nash pour une famille d'hypersurfaces quasi-rationelles, Analles de faculte des sciences de Toulouse, vol.XX, 3, 2011,
   \cite{MR1395174}, \cite{MR1460599},
    \cite{MR2153965},  \cite{petrov-2006},
    \cite{MR3109732}, \cite{MR2349665}, \cite{MR2911273}, and others). An important result about the geometry of arc spaces is the theorem of Grinberg - Kazhdan - Drinfeld, describing the formal neighborhood of a closed arc in $X_{\infty}$.\\


In this article we are trying to review the known facts about, or related to, the formal neighborhoods in arc spaces. The second section represents the theorem of Grinberg - Kazhdan, originally proved over $\mathbb{C}$ (\cite{MR1779611}). Drinfeld extended it over arbitrary field $k$, giving a purely algebraic proof in his paper (Drinfeld - On the Grinberg - Kazhdan formal arc theorem, ArXiv: math/0203263v1, see also C.Bruschek, H.Hauser - Arcs, cords and felts - six instances of the linearization principle, American Journal of Mathematics, 132, 2010). Our purpose is to give a complete proof with details and explanations, often missing in the original paper, with an example, generalizing slightly the one given by Drinfeld. In the third section we discuss two possible approaches to the relative case of the theorem.\\

\textbf{Acknowledgements}. I would like to thank cordially to Javier de Bobadilla, Shihoko Ishii and Mircea Mustata for the stimulating questions, constructive critics and helpful remarks, and to Roi Docampo for the numerous crucial suggestions and plausible comments. My special thanks are to Mark Spivakovsky for reading thoroughly the draft, and for the invaluable constructive critics, which have been vital for this article.\\

\textbf{Notations}:
1) $X$ is a scheme of finite type over a field $k$, with singular locus $X_{sing}$;\\
2) $X_{\infty}$ is the space of arcs of $X$, parameterizing all morphisms\\
$\spec (K[[t]]) \rightarrow X$, $K/k$ being an extension;\\
3) $FN_{\gamma_0}$ is the formal neighbourhood of an arc ${\gamma}_0 \in X_{\infty}$, i.e. $\spf (\widehat{{\mathcal{O}}_{X_{\infty}, {\gamma}_0}})$, the completion being taken with respect to the maximal ideal $m_{\gamma_0}$.\\
4) As a topological space, $\spec (k[[t]]) = \{0, \eta\}$, with $0$ the closed point, and $\eta$ the generic one.\\
5) $D^{\infty}$ is the formal scheme defined by the ring $k[[x_i, i \in {\mathbb{N}}]]$ of all formal expressions $\sum_{i=0}^{\infty} f_i$, each $f_i$ a homogeneous polynomial of degree $i$ in the variables $x_j$. This ring is the completion of $k[x_i, i \in \mathbb{N}]$ by the $(x_i, i \in \mathbb{N})$-adic topology, viewed as an object in the category of topological rings and continuous homomorphisms (from now on we write $(x)$ for $(x_i, i \in {\mathbb{N}})$). So, it is complete local ring with maximal ideal $m$, with the topology of projective limit, and a base for the open neighborhoods of 0 given by the decreasing sequence of ideals $I_n$, consisting of all power series with $f_i = 0, i \leq n$. The defined linear topology is weaker than the $m$-adic one, and though the ring $k[[x_i, i \in {\mathbb{N}}]]$ is complete in it, it is not complete in the $m$-adic topology. Moreover, $(x) k[[x]]  \subsetneq m$, due to the fact that $k[[x]]$ is a completion of a non Noetherian ring. \\

\section{Grinberg - Kazhdan - Drinfeld theorem}

In this section is given a complete proof of the main theorem over arbitrary field $k$. For it we need the following
\begin{defn} A test ring $A$ is a local $k$-algebra with a unit (not necessary Noetherian), having residue field $k$, and nilpotent maximal ideal $m$ (i.e. $m^n = 0$ for some $n \in \mathbb{N}$). Let $\mathit{Testrings}_k$ be the category of test rings with local homomorphisms.
\end{defn}
A typical example is $k[t]/(t^n)$.

\begin{thm} (Grinberg-Kazhdan-Drinfeld) If $\gamma_0 \in X_{\infty} \backslash (X_{sing})_{\infty}$ is a $k$-arc, there exists a scheme of finite type $Y$ and a point $y \in Y(k)$ with formal neighbourhood $FN_y = \spf \widehat{\mathcal{O}}_{Y, y}$, so that the formal neighbourhood $FN_{\gamma_0}$ is isomorphic to $\spf{\widehat{\mathcal{O}}_{Y, y}}[[{z_i}, i \in \mathbf{N}]]$, the completion of ${\widehat{\mathcal{O}}_{Y, y}}[{z_i}, i \in \mathbf{N}]$ with respect to the ideal $({z_i}, i \in \mathbf{N})$. That is, $FN_{\gamma_0} \simeq FN_y\hat{\times} D^{\infty}$.\\
\end{thm}

We start with some remarks:\\

1) The ring $\widehat{\mathcal{O}}_{Y, y}[[{z_i}, i \in \mathbf{N}]]$ consists of all power series of the form $\sum g_n$, with $g_n$ homogeneous polynomials in variables $z_i$ of degree $n$ with coefficients in $\widehat{\mathcal{O}}_{Y, y}$.\\

2) If $\gamma_0 : \spec k[[t]] \rightarrow X$, and $X'$ is the closure of the irreducible component of $X_{reg}$ containing $\gamma(\eta)$, then $X_{\infty,{\gamma}_0} = X'_{\infty,{\gamma}_0}$, so we can assume $X$ to be reduced and irreducible. Also, because the claim is local we could take $X$ to be affine as well. But the scheme $Y$ does not need to be neither reduced nor irreducible in general (see Example 2.5).\\

3) If $S$ is arbitrary (not necessarily Noetherian) scheme over $k$, and $p \in S(k)$, the formal neighbourhood $FN_p$ could be defined by its ``restricted'' functor of points\\
$h_S : Testrings_k \rightarrow Sets$,\\
$A \mapsto FN_p(A)$, the set consisting of all morphisms $\spec A \rightarrow FN_p$.\\
The reason that the category $Testrings_k$ is enough to define the functor of points is that the complete local ring $\widehat{\mathcal{O}}_{S, P}$ is a projective limit of test rings. That is, for any complete local ring $R$,\\ $Hom(R, \widehat{\mathcal{O}}_{S, P}) = lim_{\leftarrow n}Hom(R, \widehat{\mathcal{O}}_{S, P}/m_{P}^n)$, so one could apply Yoneda's lemma in the category of complete separated local $k$-algebras with linear topology weaker than the $m$-adic topology, to the category of functors $Testrings_k \rightarrow Sets$. In particular, for a test ring $(A, m_A)$ the $A$-points of $FN_{\gamma_0}$ correspond to the $A[[t]]$-points on $X$, whose reduction modulo $m_A$ is $\gamma_0$, i.e. to the $A$-deformations of $\gamma_0$.\\

Proof of Thm 2.2. Let $X$ be reduced, irreducible of dimension $n$, embedded in ${\mathbf{A}}^N$.\\
{\it{Claim 1}}. When working with local properties of $\pi^{-1}(X_{sing})$,\\
$\pi: X_{\infty} \rightarrow X$ being the canonical projection, we could take without loss of generality $X$ to be a locally complete intersection, which may be reducible. Indeed, let $r = \cod X$, and let the ideal of $X$ be $I_X = \{f_1,\dots f_s\}$. Define $F_i = \sum{a_{ij}f_{j}}$, $i=1, \dots s$, with $a_{ij}$ being general coefficients in $k$, and let $M \subset \mathbf{A^N}$ be the ideal $I_M = (F_1, \dots F_r)$ defined by the first $r$ of them. Then the following hold:\\
1) any irreducible component of $M$ has dimension $n$, so $M$ is a complete intersection scheme;\\
2) $X \hookrightarrow M$ is a closed subscheme, and $X$ and $M$ coincide at the generic point of $X$, that is, on an open nonempty subset;\\
3) there is some $r$-minor of the Jacobian matrix of $M$ not vanishing at $\eta_X$;\\
4) $X_{sing} \subset M_{sing}$.\\
Clearly, the same properties will hold for any choice of $r$ among the $F_i$'s.\\
{\it{Claim 2}}. There exists closed affine complete intersection scheme of finite type $X' \supset X$ of the same dimension such that $\im (\gamma_0)$ is not contained in $\overline{X'\setminus X}$.\\
Indeed, take $L$ to be the index set of all the $r$-tuples $(i_1, \dots i_r)$ of distinct integers with $i_j \in\{1, \dots s\}$, and let $M_l$, for $l \in L$, be the corresponding complete intersection scheme. If there is no such an $X'$ as claimed, for all $M_l$ we would have $Im (\gamma_0) \subset \overline{M_l \setminus X}$, thus $Im (\gamma_0)$ is contained in their intersection. But $\overline{M_l \setminus X}\cap X \subset \sing(M_l)$, and $\bigcap_l \sing(M_l) = X_{sing}$, contradicting the choice of $\gamma_0$.\\
For such an $X'$ we have $FN^{X}_{\gamma_0} = FN^{X'}_{\gamma_0}$, because $\gamma_0(\eta) \in X_{reg}$, so without loss of generality we can replace $X$ by $X'$.\\

Now $X$ is supposed to be complete intersection affine variety, contained in $\spec k[x_1, \dots x_n, y_1, \dots y_r]$, and defined by equations\\
$p_i = 0, i=1, \dots r$. Also, $\gamma_0(t) = (x^o(t), y^o(t)) =\\
(x^o_1(t), \dots, x^o_n(t), y^o_1(t), \dots, y^o_r(t))$ is not contained for any $t$ in $X_{sing}$, defined by $det(\partial p_i / \partial y_j) = 0$.\\
For $A$ a test ring, let $\gamma = (x(t), y(t))$, with $x(t) \in A[[t]]^n,\\ y(t) \in A[[t]]^r$, be an $A$-deformation of $\gamma_0$, i.e. its reduction modulo $m \subset A$ is equal to $\gamma_0$. Because $Im (\gamma_0) \not\subset Z(\det(\partial p_i / \partial y_j))$, not all coefficients of the power series $\det(\partial p_i / \partial y_j))(\gamma_0)$ are 0. So we may apply:
\begin{lem} (Weierstrass preparation theorem) Let $(R, m)$ be complete local separated ring with respect to a linear topology, which is weaker than the $m$-adic topology, $f= {\sum}_i c_it^i \in R[[t]]$, with not all $c_i \in m$. If $d$ is the first index such that $c_i \notin m$, then we have unique representation $f = q.u$, for some monic polynomial\\ $q = t^d + \sum_{0 \leq l<d} a_lt^l  \in R[t]$ of degree $d$, with $a_l \in m$ for all $l$, and $u \in R[[t]]$ invertible.
\end{lem}
Thus, $det(\partial p_i / \partial y_j))(x(t), y(t)) = q(t)u(t)$, for some invertible $u(t) \in A[[t]]$, and $q(t) \in A[t]$ a monic polynomial of degree $d$ whose reduction modulo $m$ equals $t^d$. The degree $d$ depends on $\gamma_0$ only, not on the choice of its deformation $\gamma$. We may assume $d \geq 1$ because, if $d = 0$ we can eliminate $y$, and the claim holds.\\
The idea of the proof is to consider $q$ as a new unknown variable. Then all $A$-deformations of $\gamma_0$ are in one-to-one correspondence with the solutions of the following system of equations with unknowns $q \in A[t], x \in A[[t]]^n, y \in A[[t]]^r$:

(I)\\
$\det(\frac{\partial p_i}{\partial y_j})(x, y) = 0  \  mod \ q$,\\
$p(x, y) = (p_1, \dots p_r) = 0. $.\\

Here, if the first equation holds, $q^{-1}\det(\partial p_i / \partial y_j)$ is invertible because it is invertible modulo $m$ and $m$ is nilpotent.\\
Next, for any fixed $e \in \mathbf{N}$, consider the following system with unknowns $q \in A[t], x \in A[[t]]^n, \bar{y} \in A[[t]]^r/(q^e)$, such that $q$ is monic polynomial of degree $d$, $q = t^d \ mod \ m$, $x = x^o \ mod \ m$, $\bar{y} \ mod \ m = y^o \ mod \ t^{ed}$:

(II)\\
$\det(\frac{\partial p_i}{\partial y_j})(x, \bar{y}) = 0  \  mod \ q$,\\
$p(x, \bar{y}) \in Im(q^e(\frac{\partial p_i}{\partial y_j})(x, \bar{y}): \frac{A[[t]]^r}{q A[[t]]^r} \rightarrow \frac{q^e A[[t]]^r}{q^{e+1}A[[t]]^r})$.\\

The second condition makes sense, if one takes the Taylor expansion of $p(x, y)$ and noting that $p(x, \bar{y})$ is well defined modulo $Im(q^e(\frac{\partial p_i}{\partial y_j})(x, \bar{y}))$. Moreover, it is equivalent to the equation $\hat{C} p(x,y) = 0 \ mod \ q^{e+1}$, where $\hat{C}$ is the adjoint matrix to $C = (\frac{\partial p_i}{\partial y_j}(x, y))$ (i.e. $C\hat{C} = \det(C).I_r$), with $y = y(t) \in A[[t]]^r$ in the pre-image of $\bar{y}$. Indeed, if $p(x, \bar{y}) = q^eC(x, \bar{y}).v$ for some $v \in A[[t]]^r/(q)$, then $\hat{C}p(x, \bar{y}) = q^{e+1}v$. Conversely, if $\hat{C}p(x, \bar{y}) =0 \mod {q^{e+1}}$, there is some $w$ such that $\hat{C}p(x, \bar{y}) = q^{e+1}w$, so $C\hat{C}p(x, \bar{y}) = Cq^{e+1}w$. Then for some invertible $u$, $p(x, \bar{y}) = Cq^{e}uw$, that is $p(x, \bar{y}) \in \im q^{e}C$.\\
Furthermore, for any fixed $e \in \mathbf{N}$ the last condition in (II) is equivalent to the following equations:

$p(x, \bar{y}) = 0 \ mod \ q^e$;

$\hat{B}p(x, \bar{y}) = 0 \ mod \ q^{e+1}$, where $B = (\frac{\partial p_i}{\partial y_j}(x, \bar{y}))$.

Both come from the second condition in (II), and the second equation makes sense once the first one holds. In this way the system (II) is equivalent to the following system which does not need any choice of $y \in A[[t]]^r$ such that $y \mod q^e = \bar{y}$, and $x(t)$ is relevant up to $\bar{x} = x \ mod \ q^{e+1}$:\\

(III)\\
$\det(\frac{\partial p_i}{\partial y_j})(x, \bar{y}) = 0  \  mod \ q$,\\
$p(x, \bar{y}) = 0 \ mod \ q^e$;\\
$\hat{B}p(x, \bar{y}) = 0 \ mod \ q^{e+1}$, where $B = (\frac{\partial p_i}{\partial y_j}(x, \bar{y}))$.\\

\begin{lem} For any $e$, the natural map from the set of solutions over $A$ of system (I) to the set of solutions over $A$ of system (II) is bijective.
\end{lem}
Proof. Let $c \in \mathbf{N}$ be the minimal number such that $m^c = 0$; we will prove the lemma by induction on $c$. If $c=1$, i.e. $m = 0, A = k$, this holds because the both systems have one solution.  Let $c \geq 2$, and suppose the claim holds for $c-1$, that is, there exists $\tilde{y} \in A[[t]]^r$ with $\tilde{y} \mod q^{e} = \bar{y}$ and $p(x(t, \tilde{y}(t)) \in m^{c-1}[[t]]^r$. This is the second equation of system I)taken over $A/m^{c-1}$. Such an $\tilde{y}$ is unique modulo $q^eA[[t]]^r \cap m^{c-1}[[t]]^r$. To prove that the map between the set of solutions of the system (I) and the set of solutions of the system (II) is bijective we have to prove it has an inverse. That is, we have to find $z(t) \in q^eA[[t]]^r \cap m^{c-1}[[t]]^r$ s.t. $p(x, \tilde{y} - z) = 0$. As before, letting $C := (\frac{\partial p}{\partial y})(x(t), \tilde{y}(t))$, we would have $p(x, \tilde{y}) - C.z = 0$, because $m^{2c-2} = 0$ by the assumption $c \geq 2$. By the first equation in (II) and the argument following the system (I), $\det(C) = q.u$ for some invertible $u(t) \in A[[t]]$. This means that $z(t)$ will be unique if it exists, because $C.z = p(x, \tilde{y})$, and $\det(C) \neq 0$. By the second equation in (II), $p(x, \tilde{y}) \in q^{e}C.A[t]^r + q^{e+1}A[t]^r$. But $C.A[t]^r \supset C.\hat{C}.A[t]^r = \det(C)I_{s}.A[t]^r = q.u.A[t]^r$, so $q^{e+1}A[t]^r = q.q^{e}A[t]^r \subset q^eC.A[t]^r$. Thus $p(x, \tilde{y}) \in q^eCA[t]^r$, that is, there is $z \in q^eA[t]^r$ such that $p(x, \tilde{y}) = C.z$. It remains to prove that $z \in m^{c-1}[t]^r$. By the induction we have $C.z = p(x, \tilde{y}) = 0 \mod m^{c-1}$, and multiplying both sides by $\hat{C}$ we have $\det(C)I_{s}.z = q.u.I_{s}.z = 0 \mod m^{c-1}$. Thus $q.z = 0 \mod m^{c-1}$, and as $q$ is monic, $z = 0 \mod m^{c-1}$ as expected.\\

We are continuing the proof of Thm 2.2. From Lemma 2.3, the $A$-deformations of $\gamma_0$ are in one-to-one correspondence with the set of solutions of the system (II), and thus, with the set of solutions of the system (III) (for any fixed $e \in \mathbf{N}$). This latter set is, in fact defined by finite number of equations in finitely many variables, because $x(t)$ could be replaced with $\bar{x} = x(t) \mod q^{e+1}$. Take $e = 1$, for example, so that $x(t) = q^2.\xi + \bar{x}$, where $\xi \in A[[t]]^n$, $\bar{x} \in A[t]^n$, $\deg(\bar{x}) < 2d$. We can consider $\bar{x}, \bar{y}, \xi, q$ as new system of unknowns, replacing $x, \bar{y}, q$. Then (II) becomes a finite system of equations over $k$ for $\bar{x}, \bar{y}, q$, with $\xi$ not involved. By the remark about the restricted functor of points above, this proves that the formal neighbourhood of $\gamma_0$ is $FN_{\gamma 0} \simeq \spf(R[[z_i, i \in {\mathbf{N}}]])$, where $R$ is a complete local Noetherian ring which defines the formal neighbourhood of a point on a scheme of finite type $y \in Y(k)$. $FN_y$ is defined by equations including the variables $\bar{x}, \bar{y}, q$ in terms of its functor of points. Thus for any $k$-algebra $S$, if $D:= (\frac{\partial p}{\partial y})(\bar{x}, \bar{y})$, let\\
$Y(S):= \{(q, \bar{x}, \bar{y}): q \in S[t], \bar{x} \in S[t]^n/(q^2), \bar{y} \in S[t]^s /(q):\\
\det(D) = 0 \mod q, p(\bar{x}, \bar{y}) = 0 \mod q, \hat{D}.p(\bar{x}, \bar{y}) = 0 \mod q^2\}$.\\
Moreover, the point $y \in Y(k)$ will correspond to $(q = t^d, \bar{x} = x^{0}(t) \mod (t^{2d}), \bar{y} = y^{0}(t) \mod (t^d)$.\\
The second factor is $D^{\infty} = \spf k[[z_i, i \in {\mathbf{N}}]]$, because $\xi \in A[[t]]^n$ has no restriction. This completes the proof.

\begin{Rem}
1) The number $d = deg(q)$ in the proof of Thm. 2.2 is equal to $\ord_{\gamma_0}(Jac_X)$, if the variables $y_j, j = 1,\dots r$ are chosen in such a way among all variables, that the $r \times r$ minor defined by them in the Jacobian matrix of $X$ has the minimal possible order among all such minors, which define locally the ideal of $X_{sing}$.\\
2) The scheme of finite type $Y$ is not unique, because one could ``enlarge'' it replacing by $Y \hat{\times} D^n$ for any $n$.
\end{Rem}

\begin{Ex} Let $\ch k = 0$ and $X: f(x_1,\dots, x_n) + x^{s}_{n+1}y = 0$ be a hypersurface in ${\mathbf{A}}^{n+2}$ with coordinates $(x_1, \dots, x_{n+1}, y)$, for $f \neq 0$ a polynomial and $s \geq 1$ an integer. Suppose $\gamma^0 := (0, \dots, 0, t, 0) \in X_{\infty}$, viewed as an $(n+2)$ -tuple of formal power series, satisfying the equation of $X$. Any $A$-deformation $\gamma$ of $\gamma_0$ is an $(n+2)$-tuple of power series $(x_1(t), \dots, x_{n+1}(t), y(t))$ satisfying\\
(1) $x_i(t) \in m[[t]], i = 1, \dots n$, $y(t) \in m[[t]]$,\\
where $m \subset A$ is the maximal ideal of the test ring $A$. By Weierstrass division theorem, any $A$-deformation of $x^0_{n+1}(t) = t$ will be of the form $x_{n+1}(t) = (t- \alpha).u(t)$ for some $\alpha \in m$ and $u(t) \in A[[t]]$ invertible. Now, given $\alpha, u(t), x_1(t), \dots x_n(t)$, there will be at most one $y(t)$, satisfying (1), and it exists iff\\
(2) $f(x_1(\alpha), \dots x_n(\alpha)) = f'(x_1(\alpha), \dots x_n(\alpha)) = \dots\\ f^{(s-1)}(x_1(\alpha), \dots x_n(\alpha)) = 0$\\
(the derivation is of power series with respect to $t$).\\
Indeed, if a solution $y(t)$ of $f(x_1(t), \dots x_{n}(t)) +x^{s}_{n+1}(t)y(t) = 0$ exists then it is unique, and $(t - \alpha)^s$ divides $f(x_1(t), \dots x_n(t))$. Conversely, if for an $\alpha \in m$, $f(\alpha) = f'(\alpha)= \dots = f^{(s-1)} (\alpha) =0$, then $y(t)$ exists because $(t- \alpha)^s$ divides $f(x_1(t), \dots x_n(t))$.\\
The system (2) defines a scheme of finite type $Y$ with $k$-point $y = (0, \dots 0)$.
\end{Ex}

\section{The relative case of Grinberg-Kazhdan-Drinfeld theorem}

Some questions about GKD theorem arise naturally. What is the geometric meaning of the scheme of finite type $Y$, appearing as first factor in it (M. Lejeune-Jalabert)? How would look like the relative case of the theorem (J. de Bobadilla)? In this section we will discuss two possible points of view, for morphisms with smooth domain, for smooth morphisms, and for closed embeddings in a smooth variety, among the others.\\

\begin{Rem} In the system (III) of the proof of Thm 2.2 we have coordinates $\bar{x}, \bar{y}, q$, and we take $e = 1$. Let $X$ be embedded in $\mathbb{A}^{n+r},\\
n = dim_P X$, $r = codim_P (X)$, and $d = deg(q) = \ord_{\gamma_0} (Jac_X)$. As the conditions of (III) give $d+rd+2rd$ equations on the $2nd+rd+d$ coefficients of $\bar{x}, \bar{y}, q$, we have the following bounds for the dimension of $Y$:\\
$d(2n+r+1) - d(1+r+2r) = d(2n-2r) \leq \dim(Y) \leq d(2n + r + 1)$. We note that the bounds hold only for a scheme $Y$ constructed in the proof of Thm.2.2. Would be interesting to understand the possible relation between any two schemes of finite type that could appear as the first factor in $FN_{\gamma_0}$.\\
 Intuitively, the more singular is $\gamma_0$, the bigger dimension has the scheme $Y$, chosen as the first factor in the formal neighborhood of $\gamma_0$. If $d= \ord_{\gamma_0}(Jac_X)=0$, $y$ could be eliminated, so $Y \simeq \spf k[[x_1, \dots x_n]]$, and the first factor could be chosen to be $\spec k$. The following example shows that the opposite does not necessary hold.
\begin{Ex}
As a particular case of Ex.2.6, let $X: xy = 0 \subset {\mathbb{A}}^2$, and let $\gamma_0 = (t, 0)$. Because $deg\ q = 1$, for the deformation $\gamma$ we have  $x(t) = u(t)(t-\alpha), y(t) = \beta$, with $u(t)$ invertible, and
$xy = 0 \ mod \ q^2$. Thus $\beta = 0$, with no restriction for $\alpha \in m_A$, the maximal ideal of a test ring $A$, and for $u$. Using the argument in Ex.2.6, we conclude that in the formal neighborhood of $FN_{\gamma}$, the scheme of finite type $Y$ could be chosen to be $\spf k[[t]]$, or just $\spec k$. That is, even when $Y$ could be chosen smooth, it could happen $\gamma(0) \in X_{sing}$.
\end{Ex}

\end{Rem}
Now suppose we have a morphism between algebraic varieties\\
$f: X \rightarrow Y$, inducing morphism $f_{\infty} \colon X_{\infty} \rightarrow Y_{\infty}$ between the arc spaces. Given $k$-arcs $\gamma \in X_{\infty} \backslash (X_{sing})_{\infty}, \delta \in Y_{\infty} \backslash (Y_{sing})_{\infty}$ with $f_{\infty}(\gamma) = \delta$, there is induced morphism $\hat{f}: FN_{\gamma} \rightarrow FN_{\delta}$ between the formal neighbourhoods. By Thm 2.2, the first one is isomorphic to $\spf \widehat{\mathcal{O}}_{U,u}[[x_i, i \in {\mathbf{N}}]]$ for a scheme of finite type $U$ and point $u \in U(k)$, and the second one is isomorphic to $\spf \widehat{\mathcal{O}}_{V,v} [[z_i, i \in {\mathbf{N}}]]$ for a scheme of finite type $V$ and point $v \in V(k)$. Let $R_1 := \widehat{\mathcal{O}}_{U,u}, R_2 := \widehat{\mathcal{O}}_{V,v}$, and let\\
$F: R_2[[z_i, i \in {\mathbf{N}}]] \rightarrow R_1[[x_i, i \in {\mathbf{N}}]]$ be the homomorphism, corresponding to
$\hat{f} :FN_u \hat{\times}D^{\infty} \rightarrow FN_v \hat{\times}D^{\infty}$.\\
.
\begin{Prop}
The morphism $\hat{f}$ defines naturally a pair of morphisms, $\Phi$ between Noetherian formal schemes, and $\Psi: D^{\infty} \rightarrow D^{\infty}$ between the second factors in $FN_{\gamma}, FN_{\delta}$.
\end{Prop}
Proof.
By restricting $F$ we get homomorphism $\phi \colon R_2 \rightarrow F(R_2) \widehat {\bigotimes}_k R_1$, which defines a morphism $\Phi$ between the corresponding formal schemes, which, composed with closed embedding, is a restriction of $\hat{f}$. The composition\\
$\psi \colon k[[z_i, i \in {\mathbf{N}}]]  \rightarrow R_2[[z_i, i \in {\mathbf{N}}]] \rightarrow R_1[[x_i, i \in {\mathbf{N}}]] \rightarrow k[[x_i, i \in {\mathbf{N}}]]$, where the first one is the natural inclusion, the second is $F$, and the last takes any power series modulo the maximal ideal $m_1 \subset R_1$, would define a morphism $\Psi: D^{\infty} \rightarrow D^{\infty}$.\\


 In this way the pair of morphisms does not give much of geometric information, and in general it is not possible to obtain back $\hat{f}$ from it. Would be preferable if $\hat{f}$ could be uniquely determined by a pair of morphisms. The first one, $\Phi$, is between $FN_u$, ``extended'' by the pre-image of $FN_v$, and $FN_v$. The second one, $\Psi$, is between the second factors $D^{\infty}$ in the formal neighborhoods. More precisely, we would like that the corresponding homomorphism splits as $F = p \circ (\phi \hat{\times} \psi)$, where\\ $\phi \colon R_2 \rightarrow F(R_2) \widehat {\bigotimes}_k R_1$ is a homomorphism between complete local Noetherian $k$-algebras, and $\psi \colon k[[z]] \rightarrow k[[y]]$. Taking\\ $p\colon F(R_2) \widehat{\bigotimes}_k R_1[[y]] \rightarrow R_1[[x]]$ to be the natural multiplication will give $\hat{f} = (\Phi \hat{\times} \Psi) \circ \tilde{p}$, where the natural closed embedding\\
 $\tilde{p}: FN_{\gamma} \rightarrow \spf (F(R_2)\widehat{\bigotimes}_k R_1) \hat{\times} D^{\infty} \simeq \spf F(R_2) \widehat{\times}_k FN_{\gamma}$ is composed with
 $\Phi \hat{\times} \Psi \colon \spf F(R_2) \widehat{\times}_k FN_{\gamma} \rightarrow FN_{\delta}$.\\
 If this holds for any $\gamma, \delta$ as above, we could view it as the first version of Thm 2.2, which we call GKD theorem for morphisms.\\


For which classes of morphisms $f$ could $\hat{f}$ be represented in this way? We will prove it for morphisms with a smooth $X$, for \'{e}tale and smooth morphisms, and for closed embeddings in smooth variety (all conditions satisfied locally in some neighborhoods of $P = \gamma(0)$ and of $Q = \delta(0)$).\\

We start with the case of smooth $X$.
\begin{Prop}
 For any $f \colon X \rightarrow Y$ with $P = \gamma (0) \in X_{reg}$, the GKD theorem for morphisms holds for $f$. In particular, it holds for any resolution of singularities of $Y$.\\
\end{Prop}
Proof. This time $F: R_2[[z_i, i \in \mathbb{N}]] \rightarrow k[[x_j, j \in \mathbb{N}]]$. We define $\phi: R_2 \rightarrow F(R_2)$, and we let $\psi =  F \circ i$, where\\
$i: k[[z_i, i \in \mathbb{N}]] \hookrightarrow R_2[[z_i, i \in \mathbb{N}]]$. Finally, the closed embedding is defined by the natural homomorphism $p: F(R_2) \widehat{\bigotimes}_k k[[x_j, j \in {\mathbf{N}}]] \rightarrow k[[x_i, i \in {\mathbf{N}}]]$. Obviously, $F$ is fully determined by $(\phi, \psi, p)$.\\

Another class of morphisms for which this first version of GKD theorem holds is that one of smooth morphisms. For the proof we need first the case of an \'{e}tale morphism.

\begin{Prop}
If $f$ is an \'{e}tale morphism the GKD for morphisms holds for $f$.
\end{Prop}
Proof. If $f$ is \'{e}tale, the commutative diagram formed by $f: X \rightarrow Y$ and $f_{\infty}: X_{\infty} \rightarrow Y_{\infty}$ with the two canonical projections\\
$\pi_X \colon X_{\infty} \rightarrow X, \pi_Y \colon Y_{\infty} \rightarrow Y$, is Cartesian. Thus\\ $f_{\infty} \colon X_{\infty} \simeq \spec \mathcal{O}_{X, P} \times Y_{\infty} \rightarrow Y_{\infty}$ is the projection, giving isomorphism
$h \colon R_1[[x_i, i \in \mathbb{N}]] \rightarrow  \widehat{\mathcal{O}}_{X,P} \hat{\otimes}R_2[[z_i, i \in {\mathbf{N}}]]$. So the homomorphism $F$ is the composition of $h^{-1}$ with the natural inclusion\\
$i \colon R_2[[z_i, i \in {\mathbf{N}}]] \rightarrow \widehat{\mathcal{O}}_{X,P} \hat{\otimes}R_2[[z_i, i \in {\mathbf{N}}]]$.  Define $\phi$ to be\\
$R_2 \hookrightarrow \widehat{\mathcal{O}}_{X,P} \hat{\otimes}R_2$, $\psi = id_{k[[z]]}$, and let $p = h^{-1}$. It is clear that they determine $F$ completely.\\

Then, we get the case of smooth morphism.

\begin{thm}
The GKD for morphisms holds when $f$ is smooth.
\end{thm}
Proof. In an open neighborhood $P \in U$, $f$ could be represented as a composition of \'{e}tale $g \colon U \rightarrow Y \times {\mathbf{A}}^n$, with the projection\\
$pr_1 \colon  Y \times {\mathbf{A}}^n \rightarrow Y$. Now we need the following:\\
\begin{lem}
  For any smooth $V$, and $pr_1 \colon X \times V \rightarrow X$, if\\
  $ \gamma \in (X \times V)_{\infty} \backslash  (\sing(X \times V))_{\infty}, \delta \in X_{\infty} \backslash (X_{sing})_{\infty}$ are closed arcs such that $pr_{1, \infty}(\gamma) = \delta$, then GKD for $pr_1$ holds.
\end{lem}
Proof. We have $(X \times_{k} V)_{\infty} = X_{\infty} \times V_{\infty}$, so $\gamma = (\delta, \eta)$, and $FN_{\gamma} = FN_{\delta} \hat{\times} FN_{\eta}$. By Thm. 2.2,
$FN_{\delta} = \spf R[[x_i, i \in {\mathbf{N}}]]$ for a complete local Noetherian ring $R$, thus $FN_{\gamma} = \spf (R[[x_i, i \in {\mathbf{N}}]] \hat{\otimes} k[[z_j, j \in \mathbb{N}]])$. Define $\phi \colon R \rightarrow R \hat{\otimes} R,
\phi(q) = q \otimes 1$, take $\psi$ to be the inclusion\\
$k[[x_i, i \in {\mathbf{N}}]] \hookrightarrow k[[x_i, z_j, i, j \in {\mathbf{N}}]]$, and take $p \colon  R \hat{\otimes} R\rightarrow R$ to be the natural homomorphism.\\

Continuing with the proof of Thm. 3.6., we apply Prop. 3.5 and Lem.3.7 to $F$. By taking $V = {\mathbf{A}}^n$, we get a composition of homomorphisms\\
$R_2[[u_i, i \in {\mathbf{N}}]] \rightarrow R_2[[u_i, v_j, i, j \in {\mathbf{N}}]] \rightarrow {\mathcal{O}}_{U,P} \hat{\bigotimes} R_2[[u_i, v_j, i, j \in {\mathbf{N}}]]$. Define $\phi \colon R_2 \rightarrow {R_2 \hat{\bigotimes} \mathcal{O}}_{U,P} \hat{\bigotimes} R_2, \phi(q) = q \otimes 1 \otimes 1$, and let\\
$\psi \colon k[[u_i, i \in {\mathbf{N}}]] \rightarrow k[[u_i, v_j, i, j \in {\mathbf{N}}]]$ be the natural inclusion. Taking $p \colon R_2 \hat{\bigotimes} {\mathcal{O}}_{U,P} \hat{\bigotimes} R_2 \rightarrow {\mathcal{O}}_{U,P} \hat{\bigotimes} R_2$ to be the natural homomorphism will complete the proof.\\

\begin{thm}
For any closed embedding $f \colon X \hookrightarrow {\mathbb{A}}^n$ the GKD for morphisms holds.
\end{thm}
Proof. The induced morphism on arc spaces $f_{\infty} \colon X_{\infty}\hookrightarrow {\mathbb{A}}^n_{\infty}$ defines for an arc $\gamma$ closed embedding
$\hat{f} \colon FN_{\gamma}^X \hookrightarrow FN_{\gamma}^{{\mathbb{A}}^n}$. Following the proof of Thm 2.2 and taking $e = 1$, the system of equations (III) defines scheme of finite type $(U, u)$ and the first factor $FN_u = \spf R$ of $FN_{\gamma}^X$. The corresponding to $\hat{f}$ surjective homomorphism is\\ $F \colon k[[z]] \rightarrow R[[x]]$. To define $R$ one needs
$N = d(2n+r+1)$ coordinates $z_i$ in $FN_{\gamma}^{{\mathbb{A}}^n} \simeq \spf k[[z_1, \dots z_N]][[Z]]$, where $Z = (Z_j, j \in \mathbb{N})$. So, there is natural surjective homomorphism between complete Noetherian local rings
$k[[z_1, \dots z_N]] \rightarrow R $, and we define $\phi: k[[z_1, \dots z_N]] \rightarrow R \widehat{\otimes} R$. The rest of the coordinates are the coefficients of power series, which are not under any restriction (see the end of the proof of Thm 2.2, p. 6). Let $\psi \colon k[[Z]] \rightarrow k[[x]]$ be the natural surjection, corresponding to the restriction of $\hat{f}$ on $D^{\infty}$ as second factor in $FN_{\gamma}^X$, and take $p\colon R \widehat{\otimes} R[[x]] \rightarrow R[[x]]$ to be the natural homomorphism.\\
\begin{Cor}
The same holds for any smooth variety $V$, and\\ $f \colon X \hookrightarrow V$ closed embedding.
\end{Cor}
 Proof. By the proof of Thm 2.2, in the system (I) the variables $y_j$ now could be eliminated, so the variables $x_i, i = 1, \dots, n$ become free. Then the ring defining $FN_{\gamma}^V$ is $k[[z]]$, and the surjective homomorphism $F \colon k[[z]] \rightarrow R[[x]]$ will define again a surjective $k[[z_1, \dots z_N]] \rightarrow R$, where $n = \dim X$.\\
\begin{Prop}
If $f_1 \colon X \rightarrow Y, f_2 \colon Y \rightarrow Z$ are morphisms, and\\ $\gamma \in X_{\infty} \setminus (X_{sing})_{\infty}, \delta \in Y_{\infty} \setminus (Y_{sing})_{\infty}, \eta \in Z_{\infty} \setminus (Z_{sing})_{\infty}$ are closed arcs, such that $f_{\infty}(\gamma) = \delta, g_{\infty}(\delta) = \eta$. Let $\hat{f_1}, \hat{f_2}$ be the induced morphisms between their formal neighborhoods, with GKD theorem for morphisms holding for $f_1$ and $\gamma, \delta$, and for $f_2$ and $\delta, \eta$. Then it holds for $f_2 \circ f_1$ and $\gamma, \eta$.
\end{Prop}
Proof. Suppose $F_1 \colon R'[[y]] \rightarrow R[[x]]$ and $F_2 \colon R''[[z]] \rightarrow R'[[y]]$ are the homomorphisms defined by $\hat{f_1}, \hat{f_2}$, for complete local Noetherian rings $R, R', R''$, and let $(\phi_i, \psi_i)$ define $F_i, i = 1, 2$, with appropriate surjections. Define $\phi$ to be the composition\\ $R'' \rightarrow F_2(R'') \hat{\otimes} R' \rightarrow F_1(F_2(R''))\hat{\otimes} F_1(R') \hat{\otimes} R$, which is a homomorphism between complete Noetherian rings,
 take\\ $\psi = \psi_1 \circ \psi_2 \colon k[[z]] \rightarrow k[[x]]$, and let $p$ be the natural surjection as before. Then $F_1 \circ F_2 = p \circ (\phi \hat{\otimes} \psi)$.\\

\begin{Rem}
1) What happens if $\im(\gamma) \subset X_{sing}$? One could stratify $X$ by $X_1 \colon = X_{reg}$, then $X_2 \colon = X_{sing} \setminus \sing(X_{sing})$, taking $X_{sing}$ as a subscheme of $X$, and so on. This sequence will terminate after finite number of steps. Applying Thm 2.2 to each $X_i$ in this sequence would represent the formal neighborhood of any closed arc in $X_{\infty}$, though taken as a point in a different subscheme of it.\\
2) An interesting question (J. de Bobadilla) would be if one could have a result similar to Thm 2.2 for (some class of) non-closed arcs on $X_{\infty}$ as well.
\end{Rem}

Another possible approach to the relative case of GKD theorem goes as follows. This time the arcs considered are not points in the space of arcs $X_{\infty}$, but in the relative space of arcs $J_{\infty}(X/Y)$ (\cite{MR2349665}). First we remind some basic definitions and properties.\\ Let $S$ be any ring, $R, T$ any $S$ algebras with $f: S \rightarrow R$ the natural homomorphism.
\begin{defn}
For $m \in \mathbb{N} \cup \{\infty\}$, a derivation of order $m$ from $R$ to $T$ is a sequence $(D_0, D_1, \dots D_m)$ (respectively, $(D_0, D_1, \dots ))$, where $D_0: R \rightarrow T$ is an $S$-algebra homomorphism, and $D_i: R \rightarrow T$ is a homomorphism of abelian groups, such that for all $i>0$\\
1) $D_i(f(a)) = 0$ for any $a \in S$;\\
2) $D_k(xy) = \sum_{i+j=k} D_i(x)D_j(y)$\\
The set of all derivations of order $m$ is $Der^m_{S}(R, T)$.
\end{defn}
\begin{Ex}
 For an affine variety $X = \spec k[x]/(f)$, where $x = (x_1,\dots x_n), f = (f_1, \dots f_s)$, the arc space is\\
$X_{\infty} = \spec k[x^{(0)},\dots x^{(m)}, \dots]/(f, f', \dots, f^{(m)}, \dots)$, where\\ $x^{(m)} = (x_1^{(m)}, \dots x_n^{(m)})$. Define a derivation on the ring by\\ $D(x_i^{(j)}) = x_i^{(j+1)}$, for any $j$, and any $i = 1 \dots n$, and let $D^{(p)} = D \circ \dots \circ D$, $p$ times. Then $(D^{(0)} = id, D^{(1)}, \dots )$ is a derivation of order infinity of $k[x^{(0)},\dots x^{(m)}, \dots]/(f, f', \dots, f^{(m)}, \dots)$ into itself.
\end{Ex}
\begin{defn}
The $R$-algebra of Hasse - Schmidt derivations is\\
$HS^m_{R/S} := R[x^{(i)}, x \in R, i = 0, \dots m]/I$ (respectively, if $m = \infty$,\\
$R[x^{(i)}, x \in R, i \in \mathbb{N}]/I$), where $x^{(i)}$ are variables, $x^{(0)} = x$ for any $x \in R$, and $I$ is the ideal generated by the following relations:\\
$(x+y)^{(i)} = x^{(i)} + y^{(i)}, f(a)^{(i)} = 0, i \geq 1, (xy)^{(k)} = \sum_{i+j=k} x^{(i)}y^{(j)},\\
x,y \in R, a \in S, i=1, \dots m$, for any $k \leq m$. The universal derivation $R \rightarrow HS^m_{R/S}$ is defined to be $(d_{0}, \dots d_m)$, where $d_i(x) = x^{(i)} \mod I$.
\end{defn}
Then it is easy to see that $HS^{\infty}_{R/S} = lim_{\rightarrow m} HS^{m}_{R/S}$.
\begin{Prop}
i) $HS^m_{R/S}$ with the universal derivation represents the functor $Der^m_{S}(R, .): \mathcal{A}lg_S \rightarrow \mathcal{S}ets$.\\
ii) There is a natural bijection\\
${Der^{m}_{S}(R, T) \rightarrow {\Hom}_{S}(R, T[t]/(t^{m+1}))}$, if $m < \infty$, and\\
${Der^m_{S}(R, T) \rightarrow {\Hom}_S(R, T[[t]])}$, for $m = \infty$,\\
giving by i) natural bijection\\
${Hom_{S}(HS^m_{R/S}, T) \rightarrow {\Hom}_{S}(R, T[t]/(t^{m+1}))}$, if $m < \infty$, or\\
${Der^m_{S}(HS^m_{R/S}, T) \rightarrow {\Hom}_S(R, T[[t]])}$, for $m = \infty$.\\
\end{Prop}
This allows to formulate in case of schemes the following:
\begin{thm}
There exists a sheaf $HS^m_{X/Y}$ of $\mathcal{O}_X$ algebras, such that for any open sets $\spec A \subset Y$, and $\spec B \subset f^{-1}(\spec A)$, one has\\
$\Gamma(\spec B, HS^m_{X/Y}) \simeq HS^m_{B/A}$, and these isomorphisms are compatible with localizations.
\end{thm}
\begin{defn}
The scheme of $m$ jet differentials of $X$ over $Y$ is $J_m(X/Y) = \spec HS^m_{X/Y}$, and it represents the functor
${\mathcal{S}ch}_Y \rightarrow {\mathcal{S}ets}$, $Z \mapsto Hom_Y(Z \times \spec \mathbb{Z}[t]/(t^{m+1}), X)$ (respectively,\\ $Hom_Y(Z \times \spec \mathbb{Z}[[t]], X)$).
\end{defn}
 Suppose $f: X \rightarrow Y$ with $f$ dominant, is a flat family. To avoid a component of the fiber of $f$ to be singular on $X$, we assume\\ $\dim X_{sing} < \dim X - \dim Y$ (in particular, $X_{sing} = \emptyset$ if $\dim X = \dim Y$). The relative arc space $J_{\infty}(X/Y)$ then could be viewed as the family of arc spaces on the fibers.
We note that even having $P \in X_{reg}$, the fiber over $Q$ could be singular at $P$.\\
For $\gamma \in J_{\infty}(X/Y)$ and $\delta$ closed arcs with $f_{\infty}(\gamma) = \delta$ as above, we have $\delta = \delta_Q$, the constant arc at $Q = f(P)$. The proof of Thm 2.2 for such $\gamma$ remains valid, the coefficients being in $k(Q) \simeq k$.
\begin{defn}
Let $FN^{rel}_{\gamma} \subset FN_{\gamma}$ be the formal neighborhood of the arc $\gamma \in J_{\infty}(X/Y)$.
\end{defn}
By Thm. 2.2 there are schemes of finite type, each with a $k$-point $(W, w), (U, u)$, such that the formal neighborhoods $FN_w, FN_u$ are the first factors in $FN^{rel}_{\gamma}, FN_{\gamma}$ respectively. We note that $(W, w)$ is a scheme of finite type over $\spec {\mathcal{O}_{Y, Q}}$.
By Thm.2.1 in (\cite{MR2349665}), we get the sequence
$J_{\infty}(X/Y) \rightarrow X_{\infty} \rightarrow Y_{\infty} \rightarrow Y$, in which the first morphism is closed embedding, the second morphism is $f_{\infty}$, and the last one is $\pi_Y$.The induced morphisms on formal neighborhoods, with the first one being again closed embedding, are\\
(1) $FN^{rel}_{\gamma} \rightarrow FN_{\gamma} \simeq FN_u \hat{\times} D^{\infty}\rightarrow FN_v \hat{\times} D^{\infty} \simeq FN_{{\delta}_Q}$.
Let\\ $C_Q = \spec (\mathcal{O}_X {\otimes}_{\mathcal{O}_Y} k(Q))$ be the fiber over $Q = (c_1, \dots c_l) \in Y$, so $\gamma \in (C_Q)_{\infty}$ (with some abuse of notation, because it is, in fact the image of $\gamma$ in $X_{\infty}$). From the proof of Thm. 2.2, the system of equations defining the first factor $FN_w$ in $FN^{rel}_{\gamma}$ includes all equations in system (I) defining $FN_u$, together with\\
 $f_j(x_1(t), \dots x_n(t)) = c_j, j = 1, \dots r$, where $f = (f_1, \dots f_r)$.\\
 From L.2.3 there is a monic polynomial $q'$ in this case, as we have monic polynomial $q$ appearing in the system defining $(U, u)$. The Jacobian matrix for $C_Q$ at $P$ includes all the rows of the Jacobian matrix for $X$ at $P$, with the same variables, and we conclude that $\deg q \leq \deg q'$.\\


For $\gamma \in J_{\infty}(X/Y)$, we have $\delta = \delta_Q$, and $\im (\delta) \nsubseteq Y_{sing}$, so $Q \in Y_{reg}$, giving from the proof of Thm 2.2, $FN_{\delta} = \spf k[[z_1, \dots z_m]][[z]]$, where $m = \dim Y$. In this case the second group of the variables $y_j$ could be eliminated (see p. 5), and the ring $k[[z_1, \dots z_m]]$ defines both $FN_v$ in (1), and $FN^Y_Q$. That is, for any closed arc $\gamma \in J_{\infty}(X/Y) \setminus J_{\infty}(X_{sing}/Y)$, the first factor $FN_w$ is a scheme over $FN^Y_{Q} \simeq \spf k[[z_1, \dots z_m]]$, which does not depend on $Q$.



\begin{Ex}
Take $X: x^2 - yz = 0$ as a flat family\\
$X \rightarrow \mathbb{A}^1, (x, y, z) \mapsto x$, and let $\gamma = (0, t, 0)$. As we have seen in Ex. 2.6, the first factor in $FN_{\gamma}$ could be taken as $\spec k[t]/(t^2)$, and from Ex. 3.2, in $FN^{rel}_{\gamma}$ the first factor could be taken as $\spec k$.
\end{Ex}
Now we could formulate the second version of Thm. 2.2 for morphisms.
Let $f \colon X \rightarrow X'$ be a morphism of $Y$-schemes, and let\\ $f^{rel}_{\infty} \colon J_{\infty}(X/Y) \rightarrow J_{\infty}(X'/Y)$ be the corresponding morphism of relative arc spaces.
Suppose
$\gamma \in J_{\infty}(X/Y) \setminus J_{\infty}(X_{sing}/Y),\\ \delta \in J_{\infty}(X'/Y) \setminus J_{\infty}(X'_{sing}/Y)$ be closed arcs, with $f^{rel}_{\infty}(\gamma) = \delta$, and let $\hat{f}_{rel} \colon FN^{rel}_{\gamma} \rightarrow FN^{rel}_{\delta}$ be the induced morphism on relative formal neighborhoods, with the first factors defined by complete local Noetherian rings $R_1, R_2$, respectively. Then for the corresponding homomorphism $F_{rel}$ we are asking as before, if there is a couple of homomorphisms $\phi \colon R_2 \rightarrow F_{rel}(R_2)\hat{\otimes}R_1, \psi \colon k[[z]] \rightarrow k[[x]]$ so that with the natural surjection $p \colon F_{rel}(R_2)\hat{\otimes}R_1 \hat{\otimes} k[[x]] \rightarrow R_1[[x]]$ we would have $F_{rel} = p \circ (\phi \hat{\otimes} \psi)$.\\

 When this holds for all $\gamma, \delta$ as above, we say that the relative version of GKD theorem holds for $f$. One could ask if the corresponding equivalents of Prop. 3.4, Thm. 3.6, and Thm. 3.8 are true in this version. As usual, all conditions are taken in affine open neighborhoods of $P, Q$, respectively.
\begin{thm}
Let $f \colon X \rightarrow X'$ be a morphism of $Y$-schemes, with \\ $g \colon X \rightarrow Y, h \colon X' \rightarrow Y$, and $\dim(X_{sing}) < \dim X - \dim Y,\\ \dim(X'_{sing}) < \dim X' - \dim Y$. Suppose $\gamma \in J_{\infty}(X/Y) \setminus J_{\infty}(X_{sing}/Y),\\ \delta \in J_{\infty}(X'/Y) \setminus J_{\infty}(X'_{sing}/Y)$ are $k$-arcs such that $f^{rel}_{\infty}(\gamma) = \delta$. The relative GKD holds for $f$ if:\\
i) $g$ is smooth morphism;\\
ii) $f$ is smooth morphism;\\
iii) $f$ is closed embedding, with $h$ smooth morphism.

\end{thm}
Proof. i) We want to show that there are homomorphisms\\
$(\phi, \psi, p)$ as above, which will determine uniquely the homomorphism $F_{rel}$.
If the first factor in $FN_{\gamma}$ is $\spf R$, and the first factor in $FN_{\delta}$ is $\spf R'$, we have that $FN^{rel}_{\gamma} \hookrightarrow FN_{\gamma}, FN^{rel}_{\delta} \hookrightarrow FN_{\delta}$, taken with $\hat{f}^{rel}, \hat{f}$ form a commutative diagram. The corresponding diagram of homomorphisms of rings contains then surjective homomorphism\\ $H \colon R'[[z]] \rightarrow R_2[[v]]$ and
$G \colon R[[x]] \rightarrow R_1[[u]]$. Here $R, R_1$ are the rings defining the first factors in the formal and relative formal neighborhood of $\gamma$, and $R', R_2$ are the corresponding rings for $\delta$.
As $g$ is smooth, its fibers are geometrically smooth, $f$ is a morphism over $Y$, so we could apply Prop. 3.4.

ii) In this case, we have that $f$ restricted to $C_Q$ is smooth morphism, and $\hat{f}_{rel}$ is defined by it. Then we obtain the claim by Thm.3.6.

iii) Now the restriction of $f$ on the fiber $C_Q$ is closed embedding in a smooth variety, namely, the fiber $D_Q$ of $h$, and by Cor. 3.9 the claim holds.\\

\begin{Rem}
It is easy to see that Prop. 3.10 also remains valid in this second version of GKD theorem.
\end{Rem}

\bibliographystyle{alpha}
\bibliography{peter}

\emph{Peter Petrov, Escola de Engenharia Industrial de Volta Redonda, Universidade Federal Fluminense,
Volta Redonda, RJ, 27225-125, Brazil}\\

{\it{E-mail address}}: pk5rov@gmail.com

\end{document}